\documentclass[11pt,a4paper]{article}

\usepackage{url,float,amsmath,amssymb,latexsym,pstricks,mathrsfs,comment,amsthm,graphicx,tikz,enumerate,rotating,
pgffor,cite,wrapfig,multicol,float,enumitem,dsfont,bbm,bibspacing}

\newcommand{\nc}{\newcommand}
\nc{\rnc}{\renewcommand}
\setlist[enumerate]{leftmargin=*}

\usepackage{geometry}  \rnc{\ss}{\smallskip} \nc{\ms}{\medskip} \nc{\bs}{\bigskip}

\begin{document}

\hyphenation{mon-oid mon-oids homo-morphism anti-homo-mor-phism}

\nc{\Sym}{\operatorname{Sym}}
\nc{\Alt}{\operatorname{Alt}}
\nc{\Sq}{\operatorname{Sq}}
\nc{\Atom}{A}
\nc{\Aut}{\operatorname{Aut}} \nc{\AntiAut}{\operatorname{Aut}^-} \nc{\AutAntiAut}{\operatorname{Aut}^\pm}
\nc{\End}{\operatorname{End}}
\nc{\Z}{\mathbb Z}
\nc{\comm}{\leftrightarrow}
\nc{\MG}{M_\Ga}
\nc{\FCS}[1]{X_{\comm}^+}
\nc{\FCM}[1]{X_{\comm}^*}
\nc{\GS}[1]{S_{#1}}
\nc{\GG}[1]{G_{#1}}
\nc{\xb}{\overline{x}}
\nc{\yb}{\overline{y}}
\nc{\zb}{\overline{z}}
\nc{\G}{\mathscr G}
\nc{\M}{\mathcal M}
\nc{\GL}{\operatorname{GL}}
\nc{\PGL}{\operatorname{PGL}}
\nc{\PSL}{\operatorname{PSL}}
\nc{\SL}{\operatorname{SL}}
\nc{\bbC}{\mathbb C}
\nc{\MnC}{\M_n(\bbC)}
\nc{\GLnC}{\GL(n,\bbC)}
\nc{\SLnC}{\SL(n,\bbC)}
\nc{\PGLnC}{\PGL(n,\bbC)}
\nc{\PSLnC}{\PSL(n,\bbC)}
\nc{\tr}{^{\operatorname{T}}}
\rnc{\iff}{\ \Leftrightarrow \ }
\nc{\1}{\id}
\nc{\Sone}{S^1}
\nc{\MGIJP}{\mathcal M^0(G;I,J;P)}
\nc{\MGIJQ}{\mathcal M^0(G;I,J;Q)}
\nc{\an}{\aleph_0}
\nc{\one}{\mathbbm{1}}
\nc{\KG}{K_G}
\nc{\KAutS}{K_{\Aut(S)}}
\nc{\DS}{D_S}
\nc{\der}[1]{[#1,#1]}
\nc{\KSymX}{K_{\Sym(X)}}
\nc{\DTX}{D_{\T_X}}

\nc{\itemit}[1]{\item[\emph{(#1)}]}
\nc{\E}{\mathcal E}
\nc{\TX}{\T(X)}
\nc{\TXP}{\T(X,\P)}
\nc{\EX}{\E(X)}
\nc{\EXP}{\E(X,\P)}
\nc{\SX}{\S(X)}
\nc{\SXP}{\S(X,\P)}
\nc{\Sing}{\E}
\nc{\idrank}{\operatorname{idrank}}
\nc{\SingXP}{\Sing(X,\P)}
\nc{\De}{\Delta}
\nc{\sgp}{\operatorname{sgp}}
\nc{\mon}{\operatorname{mon}}
\nc{\Dn}{\mathcal D_n}
\nc{\Dm}{\mathcal D_m}

\nc{\keru}{\operatorname{ker}^\wedge} \nc{\kerl}{\operatorname{ker}_\vee}

\nc{\C}{\mathscr C}
\rnc{\O}{\operatorname{O}}

\nc{\ds}{\displaystyle}

\rnc{\H}{\mathscr H}
\rnc{\L}{\mathscr L}
\nc{\R}{\mathscr R}
\nc{\D}{\mathscr D}
\nc{\J}{\mathcal J}

\nc{\wb}{\overline{w}}
\nc{\sm}{\setminus}
\nc{\I}{\mathcal I}

\nc{\COMMA}{,\quad}
\rnc{\S}{\mathcal S}

\nc{\T}{\mathcal T} 
\nc{\A}{\mathscr A} 
\nc{\B}{\mathcal B} 
\rnc{\P}{\mathcal P} 
\nc{\K}{\mathcal K}
\nc{\PB}{\mathcal{PB}} 
\nc{\rank}{\operatorname{rank}}

\nc{\sub}{\subseteq}
\nc{\la}{\langle}
\nc{\ra}{\rangle}
\nc{\mt}{\mapsto}
\nc{\im}{\mathrm{im}}
\nc{\id}{\mathrm{id}}
\nc{\bn}{\mathbf{n}}
\nc{\ba}{\mathbf{a}}
\nc{\bl}{\mathbf{l}}
\nc{\bm}{\mathbf{m}}
\nc{\bk}{\mathbf{k}}
\nc{\br}{\mathbf{r}}
\nc{\ve}{\varepsilon}
\nc{\al}{\alpha}
\nc{\be}{\beta}
\nc{\ga}{\gamma}
\nc{\Ga}{\Gamma}
\nc{\de}{\delta}
\nc{\ka}{\kappa}
\nc{\lam}{\lambda}
\nc{\Lam}{\Lambda}
\nc{\si}{\sigma}
\nc{\Si}{\Sigma}
\nc{\oijn}{1\leq i<j\leq n}
\nc{\oijm}{1\leq i<j\leq m}

\nc{\AND}{\qquad\text{and}\qquad}

\nc{\bit}{\vspace{-3 truemm}\begin{itemize}}
\nc{\bmc}{\vspace{-3 truemm}\begin{multicols}}
\nc{\emc}{\end{multicols}\vspace{-3 truemm}}
\nc{\eit}{\end{itemize}\vspace{-3 truemm}}
\nc{\ben}{\vspace{-3 truemm}\begin{enumerate}}
\nc{\een}{\end{enumerate}\vspace{-3 truemm}}
\nc{\eitres}{\end{itemize}}

\nc{\set}[2]{\{ {#1} : {#2} \}} 
\nc{\bigset}[2]{\big\{ {#1}: {#2} \big\}} 
\nc{\Bigset}[2]{\Big\{ \,{#1}\, \,\Big|\, \,{#2}\, \Big\}}

\nc{\pres}[2]{\la {#1} \,|\, {#2} \ra}
\nc{\bigpres}[2]{\big\la {#1} \,\big|\, {#2} \big\ra}
\nc{\Bigpres}[2]{\Big\la \,{#1}\, \,\Big|\, \,{#2}\, \Big\ra}
\nc{\Biggpres}[2]{\Bigg\la {#1} \,\Bigg|\, {#2} \Bigg\ra}

\nc{\pf}{\noindent{\bf Proof.}  }
\nc{\epf}{\hfill$\Box$\bigskip}
\nc{\epfres}{\hfill$\Box$}
\nc{\pfnb}{\pf}
\nc{\epfnb}{\bigskip}
\nc{\pfthm}[1]{\bigskip \noindent{\bf Proof of Theorem \ref{#1}}\,\,  } 
\nc{\pfprop}[1]{\bigskip \noindent{\bf Proof of Proposition \ref{#1}}\,\,  } 
\nc{\epfreseq}{\tag*{$\Box$}}

\makeatletter
\newcommand\footnoteref[1]{\protected@xdef\@thefnmark{\ref{#1}}\@footnotemark}
\makeatother

\numberwithin{equation}{section}

\newtheorem{thm}[equation]{Theorem}
\newtheorem{lemma}[equation]{Lemma}
\newtheorem{cor}[equation]{Corollary}
\newtheorem{prop}[equation]{Proposition}
\newtheorem*{conj}{Conjecture}
\newtheorem{question}{Question}
\newtheorem*{q3}{Question 3}

\theoremstyle{definition}

\newtheorem{rem}[equation]{Remark}
\newtheorem{defn}[equation]{Definition}
\newtheorem{eg}[equation]{Example}
\newtheorem{ass}[equation]{Assumption}
\newtheorem{problem}{Problem}

\title{On groups generated by involutions of a semigroup} 
\author{
James East\\
{\footnotesize \emph{Centre for Research in Mathematics; School of Computing, Engineering and Mathematics}}\\
{\footnotesize \emph{University of Western Sydney, Locked Bag 1797, Penrith NSW 2751, Australia}}\\
{\footnotesize {\tt J.East\,@\,uws.edu.au}}\\~\\
Thomas E.~Nordahl\\
{\footnotesize \emph{Department of Psychiatry and Behavioral Sciences; Imaging Research Center}}\\
{\footnotesize \emph{University of California Davis, Davis, CA 95616, United States}}\\
{\footnotesize {\tt tenordahl\,@\,ucdavis.edu}}
}
\date{For Ed Scheiblich on the occasion of his 77th birthday.}

\maketitle

~\vspace{-12ex}
\begin{abstract}
An involution on a semigroup~$S$ (or any algebra with an underlying associative binary operation) is a function $\al:S\to S$ that satisfies $\al(xy)=\al(y)\al(x)$ and $\al(\al(x))=x$ for all $x,y\in S$.  The set $I(S)$ of all such involutions on~$S$ generates a subgroup $\mathscr C(S)=\la I(S)\ra$ of the symmetric group $\Sym(S)$ on the set $S$.  We investigate the groups $\mathscr C(S)$ for certain classes of semigroups $S$, and also consider the question of which groups are isomorphic to $\C(S)$ for a suitable semigroup $S$. 


\end{abstract}

\section{Introduction}\label{sect:intro}

Involutions are ubiquitous in many branches of mathematics, and have played a particularly significant role in algebra.  
There are algebras that have (external) involution operators defined on them \cite{HLW1997,Munn1974,Scheiblich1973,Polak2001,Arveson1976,WO1993,DNN2009,Lawson1998,Bacovsky2014,Nordahl1978,GL, Higgins2014,Jones2012}, as well as algebras generated by (internal) involutions such as the well-known Coxeter Groups, mapping class groups, special linear groups, and non-abelian finite simple groups \cite{Kassabov2003,Monden2011,Humphreys1990,GHR1976,Djokovic1967,Wonenburger1966,Atlas,Kantor1973}.   An (internal) involution in a group is an element of order $2$ (i.e., a non-identity element $a$ that satisfies $a^2=1$).   
An (external) involution on a semigroup $S$ (or any algebra with an underlying associative binary operation) is a function $\al: S\to S$ satisfying $\al(\al(x)) = x$ and $\al(xy) = \al(y)\al(x)$, for every $x,y\in S$.  
Many varieties of semigroups and algebras have involutory unary operations built into their signature, including the classes of groups, inverse semigroups \cite{Lawson1998}, cellular algebras \cite{GL}, $C^*$-algebras \cite{Exel2014book,Exel2014,Arveson1976,WO1993}, MI-groups \cite{Bacovsky2014}, and regular $*$-semigroups \cite{Nordahl1978}; the latter class models (for example) several diagram monoids \cite{FitLau2011,LauFit2006,EF2012,East2011,DEEFHHL2015,EG2015,ADV2012,ADV2012_2,Maz2002,Maltcev2007}.

Some well-known algebras have multiple involutions defineable on them.   For example, the inverse and transpose maps $A\mt A^{-1}$ and $A\mt A^{\operatorname{T}}$ both define involutions on the general linear group $\GL(n,F)$, which consists of all invertible $n\times n$ matrices over a field $F$.   The composition of these two involutions (i.e., the map $A\mt(A^{\operatorname{T}})^{-1} = (A^{-1})^{\operatorname{T}}$) is a non-inner automorphism of $\GL(n,F)$.
Commuting involutions on a semisimple algebraic group (such as the inverse and transpose maps on the special linear group $\SL(n,F)\sub\GL(n,F)$) yield a $\Z_2\times\Z_2$-grading on the associated Lie algebra \cite{Panyushev2013}.
Scheiblich \cite{Scheiblich1987} gave examples of bands (idempotent semigroups) for which two involutions give rise to non-isomorphic regular $*$-semigroups.  
Auinger \emph{et al.}~\cite{ADV2012_2} studied two different involutions on the finite partition monoids (and related diagram monoids) in the context of (inherently) non-finitely based equational theories; one of these involutions leads to a regular $*$-semigroup structure and the other does not.
Winker \emph{et al.}~\cite{WWL1981} gave examples of semigroups with anti-automorphisms but no involution (they showed the minimal size of such a semigroup is $7$, and that that there are four such semigroups of minimal size, all of which are $3$-nilpotent).  
Ciobanu \emph{et al.}~\cite{CDE2015} made crucial use of a free monoid with involution in their work on word equations in free groups.
Bacovsk\'y \cite{Bacovsky2014} investigated a class of monoids with involution that have applications in processor networks and fuzzy numbers.
Gustafson \emph{et al.}~\cite{GHR1976} showed that any square matrix of determinant $\pm1$ over any field is the product of (at most) four involutions; a square matrix is the product of two involutions if and only if it is invertible and similar to its own inverse \cite{Djokovic1967,Wonenburger1966}.  A linear bound for products of involutory integer matrices was given in \cite{Ishibashi1995}.
Everitt and Fountain studied partial mirror symmetries by investigating certain factorizable inverse monoids generated by ``partial reflections'' of a space \cite{EvF2010,EvF2013}; see also~\cite{East2013_5}.
Lusztig and Vogan investigated certain actions on Hecke algebras via involutions on Coxeter groups that permute the simple reflections \cite{LV2012,LV2014}.
Francis \emph{et al.}~made extensive use of involution-generated groups in their work on algebraic bacterial genomics \cite{Francis2014,Bhatia2014,EFG2014,EGTF2014}.
Finally, we must mention the vast body of work of authors such as Dolinka, Imaoka, Jones, Petrich, Reilly, Scheiblich and Yamada \cite{Scheiblich1980,Scheiblich1982,Reilly1979,Imaoka1981_2,Imaoka1981,Imaoka1980,Imaoka1983,Yamada1981,CD2002,Dolinka2010,Dolinka2001, Dolinka2000,Dolinka2000_2,Scheiblich1987,Petrich1985,Jones2012,Jones2014} on varieties of involution semigroups and bands; see especially \cite{Petrich1985} for a discussion of early work on this topic.

In this work, we are interested in \emph{all} the involutions on a semigroup $S$, and we study these by defining a group $\mathscr C(S)$ as follows.  First, we let $I(S)$ be the set of all involutions on $S$.  So $I(S)$ is a (possibly empty) subset of $\Sym(S)$, the symmetric group on $S$ (which consists of all permutations of $S$).  We then define $\mathscr C(S)=\la I(S)\ra$ to be the subgroup of $\Sym(S)$ generated by $I(S)$.  (We interpret $\mathscr C(S)=\la\emptyset\ra=\{\1_S\}$ to be the trivial subgroup of $\Sym(S)$ in the case that $S$ has no involutions definable on it; this occurs for left- or right-zero semigroups, for example.  Throughout, we write $\1_X$ for the identity mapping on the set $X$, or just $\1$ if the set is clear from context.)

As an example, consider the Klein $4$-group, $K=\Z_2\times\Z_2$.  If we denote by $x,y,z$ the non-identity elements of $K$, then one easily checks that any permutation of $x,y,z$ induces an automorphism of~$K$.  In particular, since $K$ is commutative, the transpositions $(x,y)$, $(x,z)$ and $(y,z)$ each induce involutions of $K$.  It quickly follows that $\C(K)=\Aut(K)$ is isomorphic to $\Sym(3)$, the symmetric group on three letters.  For examples where $\C(S)$ does not coincide with $\Aut(S)$, see Sections \ref{sect:GS} and \ref{sect:families}.  In fact, it is possible for non-isomorphic semigroups $S\not\cong T$ to have $\Aut(S)\cong\Aut(T)$ but $\C(S)\not\cong\C(T)$, so the $\C(S)$ invariant distinguishes between some semigroups that the $\Aut(S)$ invariant does not (see Sections~\ref{sect:GS} and \ref{subsect:T_X}).
The main purpose of this paper is to investigate the two most natural questions to ask: 

\ms
\begin{question}\label{qu:1}
Given a semigroup $S$, can we describe the group $\mathscr C(S)$?
\end{question}

\smallskip
\begin{question}\label{qu:2}
Given a group $G$ generated by (internal) involutions, does there exist a semigroup $S$ such that $G\cong\mathscr C(S)$?
\end{question}

We answer Question \ref{qu:1} for several families of semigroups including: free semigroups and free commutative semigroups of arbitrary rank; free groups and free abelian groups of finite rank; (finite or infinite) full transformation semigroups, symmetric and dual symmetric inverse monoids; finite partition monoids; semigroups of complex matrices; rectangular bands; and graph semigroups.  

We also utilise a construction due to Ara\'ujo \emph{et al.}~\cite{ABMN2010} to give an affirmative answer to Question \ref{qu:2}.  Specifically, we show that any involution generated group is isomorphic to $\C(S)$ for some commutative semigroup $S$.  However, none of the involutions on such a commutative semigroup $S$ are ``proper'', in the sense that they only ``reverse the operation'' of $S$ (i.e., satisfy $\al(xy)=\al(y)\al(x)$) because they preserve it.
Since many semigroups of interest (such as non-commutative groups, inverse semigroups and regular $*$-semigroups) do have proper involutions, we therefore consider the following more focussed question to be of paramount importance.

\ms
\begin{question}\label{qu:3}
Given a group $G$ generated by (internal) involutions, does there exist a semigroup $S$ with a proper involution such that $G\cong\mathscr C(S)$?
\end{question}

(As alluded to above, we say an involution is \emph{proper} if it is not a homomorphism.)  Our investigation of graph semigroups will lead us to a partial answer to this question.  Namely, we show that for any involution-generated group $G$, there exists a semigroup $S$ with a proper involution such that $\C(S)\cong G\times\Z_2$.  We also show that if $G\cong\C(S)$ for some semigroup $S$ with a proper involution, then $G\cong N\rtimes\Z_2$ for some normal subgroup $N$ of $G$; in particular, this shows that Question \ref{qu:3} has a negative answer for some involution-generated groups.

{\bf Acknowledgements.}  We thank Dr Heiko Dietrich for several helpful conversations during the preparation of this article, and also Prof Mohan Putcha for encouragement and advice.

\section{Preliminaries}\label{sect:prelim}

In this section, we fix notation, and gather together various facts that we will need in our investigations.

Let $S$ and $T$ be semigroups.  Recall that a function $\al:S\to T$ is a \emph{homomorphism} (resp., \emph{anti-homomorphism}) if $\al(xy)=\al(x)\al(y)$ (resp., $\al(xy)=\al(y)\al(x)$) for all $x,y\in S$.  An anti-homomorphism is called \emph{proper} if it is not also a homomorphism (in which case $T$ must not be commutative).  An \mbox{(anti-)homomorphism} is called an (anti-)isomorphism if it is bijective.  The inverse mapping ${\al^{-1}:T\to S}$ of an \mbox{(anti-)isomorphism} $\al:S\to T$ is also an (anti-)isomorphism.  An \mbox{(anti-)isomorphism} $S\to S$ is called an (anti-)automorphism.  

We write $\Aut(S)$ (resp., $\AntiAut(S)$) for the set of all automorphisms (resp., anti-automorphisms) of the semigroup $S$.  We write $\AutAntiAut(S)=\Aut(S)\cup\AntiAut(S)$.  Note that while $\Aut(S)$ always contains the identity mapping, 
$\AntiAut(S)$ may be empty (this is true of left- or right-zero semigroups, for example, and for the full transformation semigroup $\T_X$ when $|X|\geq2$; see Section~\ref{subsect:T_X}).  The composite of two automorphisms or two anti-automorphisms is an automorphism, while the composite of an automorphism with an anti-automorphism (in any order) is an anti-automorphism.  So 
$\AutAntiAut(S)$ is a group under composition, which we call the \emph{signed automorphism group of $S$}.  

An \emph{involution} of $S$ is an anti-automorphism $\al$ of order $2$; that is, $\al\in\AntiAut(S)$ and $\al^2=\1\not=\al$.  We write $I(S)$ for the set of all involutions of $S$, and $\C(S)=\la I(S)\ra$ for the subgroup of $\AutAntiAut(S)$ generated by $I(S)$.  Note that $\C(S)=\{\1_S\}$ if $I(S)=\emptyset$.  The proof of the next result is obvious, and is omitted.

\ms
\begin{lemma}\label{lem:commutative}
Let $S$ be a semigroup.
\bit
\itemit{i} If $S$ is commutative, then $\AntiAut(S)=\Aut(S)=\AutAntiAut(S)$.
\itemit{ii} If $S$ is not commutative, then either $\AntiAut(S)=\emptyset$ or else the map $$\Phi_\be:\Aut(S)\to\AntiAut(S):\al\mt\al\be$$ is a bijection for any $\be\in\AntiAut(S)$,    
and we have a short exact sequence:
\[
\{1\} \to \Aut(S) \to \AutAntiAut(S) \to \Z_2 \to \{1\}.
\]
\eit
That is, $\Aut(S)$ is a (normal) subgroup of index at most $2$ of $\AutAntiAut(S)$. \epfres
\end{lemma}

So if $S$ is not commutative and $\AntiAut(S)$ is non-empty, then the signed automorphism group $\AutAntiAut(S)$ is an extension of $\Z_2$ by $\Aut(S)$.  This extension splits if and only if $I(S)\not=\emptyset$.  Non-commutative semigroups with $I(S)=\emptyset\not=\AntiAut(S)$ are described in \cite{WWL1981}, where it is shown that the smallest such semigroup has size $7$.

Recall that an \emph{atom} of a semigroup $S$ (which may or may not be a monoid) is an element $a\in S$ 
that cannot be factorized as a product of non-identity elements of $S$.
Write $\Atom(S)$ for the set of all atoms of~$S$ (which may be empty).  

\ms
\begin{lemma}\label{lem:atoms}
Let $\al:S\to T$ be an (anti-)isomorphism.  Then $\al(\Atom(S))=\Atom(T)$.
\end{lemma}

\pf We just prove the statement for anti-isomorphisms.
Let $a\in \Atom(S)$.  If $\al(a)=bc$ for some non-identity elements $b,c\in T$, then $a=\al^{-1}(bc)=\al^{-1}(c)\al^{-1}(b)$, where $\al^{-1}(c),\al^{-1}(b)$ are non-identity elements of $S$, contradicting the assumption that $a\in \Atom(S)$.  (Note that $S$ is a monoid if and only if $T$ is a monoid, in which case $\al$ must map the identity of $S$ to the identity of $T$.)  It follows that $\al(a)\in \Atom(T)$, so $\al(\Atom(S))\sub \Atom(T)$.  Applying this also to the anti-isomorphism $\al^{-1}:T\to S$, we obtain $\Atom(T) = \al(\al^{-1}(\Atom(T))) \sub \al(\Atom(S)) \sub \Atom(T)$, giving $\al(\Atom(S))=\Atom(T)$. \epf

We write $\Sym(X)$ for the symmetric group on a set $X$, which consists of all permutations of $X$.  In what follows, we will frequently use the fact that a symmetric group $\Sym(X)$ is generated by (internal) involutions.  Although this is certainly well known, especially in the finite case (with finite $\Sym(X)$ being a Coxeter group of Type A), we include a short proof for convenience.  (It is also known that any element of a finite Coxeter group is either an involution or the product of two involutions \cite{Carter1972}.)

\ms
\begin{lemma}\label{lem:Sym(X)}
Let $\pi\in\Sym(X)$, where $X$ is an arbitrary set.  Then $\pi=\si\tau$ for some $\si,\tau\in\Sym(X)$ with $\si^2=\tau^2=\1$.
\end{lemma}

\pf If $\pi=\1$, then the result is clear, so suppose otherwise.  Write $\pi=\prod_{i\in I} \ga_i$, where $I$ is some indexing set, and each $\ga_i$ is a non-trivial cycle of some subset $X_i\sub X$, with the $X_i$ pairwise disjoint (and each of size at least $2$).  Note that it is possible for $I$ to be infinite (in which case $\pi$ is simply the ``formal product'' of the $\ga_i$), and for some of the $X_i$ to be (countably) infinite.  Fix some $i\in I$.  We consider three cases.

{\bf Case 1.}  If $\ga_i=(x_1,\ldots,x_{2k+1})$ for some $k\geq1$, then $\ga_i=\si_i\tau_i$, where
\begin{align*}
\si_i = (x_2,x_{2k+1})(x_3,x_{2k})\cdots(x_{k+1},x_{k+2}) &\AND \tau_i = (x_1,x_{2k+1})(x_2,x_{2k})\cdots(x_k,x_{k+2}).
\intertext{{\bf Case 2.}  If $\ga_i=(x_1,\ldots,x_{2k})$ for some $k\geq1$, then $\ga_i=\si_i\tau_i$, where}
\si_i = (x_2,x_{2k})(x_3,x_{2k-1})\cdots(x_{k},x_{k+2}) &\AND \tau_i = (x_1,x_{2k})(x_2,x_{2k-1})\cdots(x_k,x_{k+1}).
\intertext{{\bf Case 3.}  If $\ga_i=(\ldots,x_{-2},x_{-1},x_0,x_1,x_2,\ldots)$ is an infinite cycle, then $\ga_i=\si_i\tau_i$, where}
\si_i = (x_1,x_{-1})(x_2,x_{-2})(x_3,x_{-3})\cdots &\AND \tau_i = (x_0,x_{-1})(x_1,x_{-2})(x_2,x_{-3})\cdots .
\end{align*}
Finally, we see that $\pi=\si\tau$, where $\si=\prod_{i\in I}\si_i$ and $\tau=\prod_{i\in I}\tau_i$ satisfy $\si^2=\tau^2=\1$. \epf

We conclude this section with a simple group theoretic result.  As usual, if $G$ is a group, we denote by $\der{G}$ the (normal) subgroup of $G$ generated by all commutators $[g,h]=ghg^{-1}h^{-1}$ with $g,h\in G$.

\ms
\begin{lemma}\label{lem:KG}
Let $G$ be a group, and denote by $\KG$ the subgroup of $G\times G$ generated by  $\set{(g,g^{-1})}{g\in G}$.  Then
\[
\KG = \bigset{(g,h)\in G\times G}{gh\in \der{G}}.
\]
In particular, $\KG$ is an extension of $G$ by $\der{G}$.
\end{lemma}

\pf Put $\Omega=\bigset{(g,h)\in G\times G}{gh\in \der{G}}$.  First note that for any $g,h\in G$, $\KG$ contains $(g,g^{-1})(hg^{-1},gh^{-1})(h^{-1},h)=(ghg^{-1}h^{-1},1)$, where $1$ denotes the identity element of $G$.  Similarly, $(1,ghg^{-1}h^{-1})\in\KG$.  It follows that $\der{G}\times \der{G}\sub\KG$.  So if $g,h\in G$ satisfy $gh\in \der{G}$, then $(g,h)=(gh,1)(h^{-1},h)\in\KG$, showing that $\Omega\sub\KG$.  

Conversely, we show by induction that
\begin{equation}\label{eq:KG}
(g_1\cdots g_k,g_1^{-1}\cdots g_k^{-1})=(g_1,g_1^{-1})\cdots(g_k,g_k^{-1})\in\Omega \qquad\text{for all $g_1,\ldots,g_k\in G$.}
\end{equation}
If $k\leq2$, this is clear, so suppose $k\geq3$.  Then
\[
g_1\cdots g_k \cdot g_1^{-1}\cdots g_k^{-1} = \big[ g_1\cdots g_{k-2}(g_{k-1}g_k)\cdot g_1^{-1}\cdots g_{k-2}^{-1}(g_{k-1}g_k)^{-1} \big] \cdot g_{k-1}g_kg_{k-1}^{-1}g_k^{-1} .
\]
By an induction hypothesis, the bracketed term belongs to $\der{G}$, giving $g_1\cdots g_k\cdot g_1^{-1}\cdots g_k^{-1}\in\der{G}$.  This completes the proof of \eqref{eq:KG}, and shows that $\KG\sub\Omega$.
Finally, we note that the epimorphism $\KG\to G:(g,h)\mt g$ has kernel $\{1\}\times \der{G}$.  \epf

\begin{rem}
The group $K_G$ may be thought of as the kernel of the epimorphism $G\times G\to G/\der{G}$ sending the pair $(g,h)$ to the coset $\der{G}gh$.  This construction may of course be generalised by replacing $G\times G$ with the direct product of $n\geq2$ copies of $G$.  The resulting groups, denoted $\mathcal K(G,n)$, are studied in \cite{Liedtke2008} in connection to invariants for complex algebraic surfaces (among other things).
\end{rem}

\ms
\begin{eg}\label{eg:KSymX}
It is well known that
\[
\der{\Sym(X)} = \begin{cases}
\Alt(X) &\text{if $|X|<\an$}\\
\Sym(X) &\text{if $|X|\geq\an$.}
\end{cases}
\]
(In fact, Ore \cite{Ore1951} showed that any element of $\der{\Sym(X)}$ is itself a commutator, and not just a product of commutators.)  It follows that
\[
\KSymX = \begin{cases}
\set{(\si,\tau)\in\Sym(X)\times\Sym(X)}{\si\tau\in\Alt(X)} &\text{if $|X|<\an$}\\
\Sym(X)\times\Sym(X) &\text{if $|X|\geq\an$.}
\end{cases}
\]
Consider now the case in which $X$ is finite and put
\[
K=\KSymX=\set{(\si,\tau)\in\Sym(X)\times\Sym(X)}{\si\tau\in\Alt(X)}.
\]
So $\Alt(X)\times\Alt(X)\leq K\leq \Sym(X)\times\Sym(X)$.  (Here, $\leq$ denotes the subgroup relation.)  In fact, we have an internal semidirect product decomposition $K=(\Alt(X)\times\Alt(X))\rtimes\la(\pi,\pi)\ra$, where $\pi$ is any fixed odd permutation of order $2$ (such as a transposition), but we note that this is not a wreath product.  We also note that the condition $\si\tau\in\Alt(X)$ is equivalent to saying that $\si$ and $\tau$ have the same parity.  
\end{eg}

\section{Obtaining any involution-generated group as $\C(S)$}\label{sect:allG}

Ara\'ujo \emph{et al.}~\cite{ABMN2010} gave a number of constructions that realise any finite group as the automorphism group of various kinds of finite semigroups.  
Here we review one of these constructions and show how it applies to the present considerations.  

Graphs play an important role in this section and the next, so we take this opportunity to fix our notation.  Throughout, graphs are assumed to be undirected and have no loops or parallel edges.  Let $\Ga$ be a graph with vertex and edge sets $X$ and $E$, respectively.  Recall that an automorphism of $\Ga$ is a bijection $\pi:X\to X$ such that $\{x,y\}\in E$ if and only if $\{\pi(x),\pi(y)\}\in E$.  We write $\Aut(\Ga)$ for the group of automorphisms of $\Ga$.  Note that $\Aut(\Ga)$ is a subgroup of $\Sym(X)$, the symmetric group over~$X$.

Let $G$ be an arbitrary group.  It is well known that $G$ is isomorphic to the automorphism group $\Aut(\Ga)$ of a suitable graph $\Ga$.  This was originally proved for finite groups (in which case the graph may be taken to be finite) by Frucht \cite{Frucht1938}, and then extended to infinite groups by de~Groot~\cite{deGroot1959}.  Various strengthenings of this result exist, including to graphs of given chromatic number or satisfying various connectivity or regularity conditions; see for example \cite{Sabidussi1957,Izbicki1959,Frucht1949}.  But all we need to know is that $\Ga$ may be chosen so that it has at least one edge (and it is easy to see why this is possible, assuming de Groot's theorem).  We fix such a graph $\Ga$, with $G\cong\Aut(\Ga)$, and write $X$ and $E$ for the vertex and edge set of~$\Ga$, respectively.  As in \cite{ABMN2010}, we define a semigroup $S=X\cup\{Y,N\}$, where $Y$ and $N$ are new symbols that do not belong to $X$, and with multiplication defined, for $u,v\in S$, by
\[
uv=\begin{cases}
Y &\text{if $u,v\in X$ and $\{u,v\}\in E$}\\
N &\text{otherwise.}
\end{cases}
\]
Note that $S$ is commutative and $3$-nilpotent (i.e., $N$ is a zero element and $abc=N$ for any $a,b,c\in S$).  It is easy to check that: 
\bit
\item[(i)] any (graph) automorphism $\al\in\Aut(\Ga)$ extends to a (semigroup) automorphism $\hat\al\in\Aut(S)$ by further defining $\hat\al(Y)=Y$ and $\hat\al(N)=N$; and 
\item[(ii)] any $\be\in\Aut(S)$ must be of the form $\hat\al$ for some $\al\in\Aut(\Ga)$.  (The assumption that $\Ga$ has an edge is necessary to prove that $\be(Y)=Y$.)
\eit
Thus, $\Aut(S)\cong\Aut(\Ga)$ and it follows that $G\cong\Aut(S)$.  In fact, since $S$ is commutative, it follows that $\AutAntiAut(S)=\Aut(S) \cong G$.  In particular, if $G$ is an involution-generated group, then so too is $\AutAntiAut(S)\cong G$, so that $\C(S)=\AutAntiAut(S)\cong G$.  This therefore proves the following.

\ms
\begin{thm}\label{thm:any_G}
If $G$ is an arbitrary involution-generated group, then $G\cong\C(S)$ for some (commutative) semigroup $S$. \epfres
\end{thm}

Theorem \ref{thm:any_G} gives an affirmative answer to Question \ref{qu:2} from the introduction.  However, since the semigroup $S$ constructed above is commutative, $S$ has no \emph{proper} involutions.  So the above considerations do not say anything about Question \ref{qu:3}, which we restate here for convenience.  

\ms
\begin{q3}
Given a group $G$ generated by (internal) involutions, does there exist a semigroup $S$ with a proper involution such that $G\cong\mathscr C(S)$?
\end{q3}

Since many well-studied involutory semigroups (such as non-commutative groups, inverse semigroups, regular $*$-semigroups, and so on) have proper involutions, Question \ref{qu:3} does seem a natural one to investigate.  The next result shows that there are some restrictions on the kind of group $G$ that may be isomorphic to $\C(S)$ for a semigroup $S$ with a proper involution; namely, $G$ must be a split extension of $\Z_2$ by some normal subgroup.

\ms
\begin{prop}\label{prop:rtimes}
Let $S$ be a semigroup with a proper involution $\iota$.  Then we have (internal) semidirect product decompositions
\[
\AutAntiAut(S) = \Aut(S)\rtimes\la\iota\ra \AND \C(S) = (\C(S)\cap\Aut(S))\rtimes\la\iota\ra.
\]
That is, $\AutAntiAut(S)$ \emph{(resp}., $\C(S)$\emph{)} is a split extension of $\Z_2$ by $\Aut(S)$ \emph{(}resp., $\C(S)\cap\Aut(S)$\emph{)}.
\end{prop}

\pf This follows quickly from the facts that $\Aut(S)$ (resp., $\C(S)\cap\Aut(S)$) is an index $2$ subgroup of $\AutAntiAut(S)$ (resp., $\C(S)$), and that $\iota$ has order~$2$. \epf

\begin{rem}
It follows immediately that the class of groups isomorphic to $\C(S)$ for some semigroup~$S$ with a proper involution is smaller than the class of all involution-generated groups.  For example, if~$G$ is a finite non-abelian simple group, then $G$ is generated by its elements of order $2$ (since $|G|$ is even, and the subgroup generated by elements of order $2$ is a non-trivial normal subgroup of $G$).  But such a group $G$ is obviously not a split extension.  
\end{rem}

We now investigate a special case in which the semidirect product decompositions in Proposition \ref{prop:rtimes} are direct products.  For a semigroup $S$, write $J(S)=\set{\al\in\Aut(S)}{\al^2=\1}$, and put $\G(S)=\la J(S)\ra$.  (Note that $\1\in J(S)$ for all $S$.)

\ms
\begin{prop}\label{prop:times}
Let $S$ be a semigroup with a proper involution $\iota$ that commutes with all automorphisms of $S$.  Then 
\[
\Psi_\iota:J(S)\to I(S):\al\mt\al\iota
\]
is a bijection, and we have (internal) direct product decompositions
\[
\AutAntiAut(S) = \Aut(S)\times\la\iota\ra \AND \C(S) = \G(S)\times\la\iota\ra.
\] 
\end{prop}

\pf If $\al\in J(S)$, then $\al\iota\in\AntiAut(S)$ and $(\al\iota)^2=\al^2\iota^2=\1$, so that $\al\iota\in I(S)$.  Conversely, if $\be\in I(S)$, then $\be=(\be\iota)\iota$, with $\be\iota\in\Aut(S)$ and $(\be\iota)^2=\1$.  To complete the proof, it suffices to show that $\C(S)\cap\Aut(S)=\G(S)$.  If $\al\in\C(S)\cap\Aut(S)$, then $\al=(\be_1\iota)\cdots(\be_{2k}\iota)$ for some $\be_1,\ldots,\be_{2k}\in J(S)$, from which it follows that $\al=\be_1\cdots\be_{2k}\in\G(S)$.  Conversely, if $\al\in\G(S)$, then $\al = \be_1\cdots\be_l$ for some $\be_1,\ldots,\be_l\in J(S)$ where we may assume that $l$ is even (since $\1\in J(S)$), in which case, $\al=(\be_1\iota)\cdots(\be_l\iota)\in\C(S)$.~\epf

In particular, Proposition \ref{prop:times} holds for non-commutative inverse semigroups (including non-abelian groups), since we may use the inverse mapping ${\iota:S\to S:s\mt s^{-1}}$, noting that $\al(s^{-1})=\al(s)^{-1}$ for any semigroup homomorphism $\al:S\to T$ between inverse semigroups.  In the next section, we will study another class of semigroups (that are not inverse semigroups) that Proposition \ref{prop:times} applies to.  As a consequence of our investigations, we will see that for any involution-generated group $G$, there exists a semigroup $S$ with a proper involution such that $\C(S)\cong\Z_2\times G$.

\ms
\begin{rem}
In the argument used to prove Theorem \ref{thm:any_G}, above, we made use of the fact \cite{ABMN2010} that every group is isomorphic to the automorphism group $\Aut(S)$ of some commutative semigroup $S$.  It is known \cite{dVdM1958,HL1969} that some groups (such as the cyclic group of order $5$) are not isomorphic to the automorphism group $\Aut(G)$ of any group $G$.  So the question of which (involution-generated) groups are isomorphic to $\C(G)$ for some group $G$ seems to be of importance.  For example, we noted in the introduction that the symmetric group of degree $3$ is isomorphic to $\C(K)$ where $K=\Z_2\times\Z_2$ is the Klein $4$-group.
\end{rem}

\section{Graph semigroups}\label{sect:GS}

In this section, we calculate $\C(S)$ in the case that $S$ belongs to the class of \emph{graph semigroups} (see below for the definition).  Part of the motivation for doing this is that the free semigroups and free commutative semigroups belong to this class.  The other motivation is that we may use the results to give a partial answer to Question 3 above.  (We note that the results of this section apply equally well to the corresponding monoids, but we just treat semigroups for convenience.)  

We begin by fixing notation for dealing with semigroup presentations.  Let $X$ be an arbitrary set, and write $X^+$ for the \emph{free semigroup} over $X$, which consists of all non-empty words over~$X$ under the semigroup operation of concatenation.  We regard $X$ as a subset of $X^+$ in the usual way (as the set of words of length $1$).  For $R\sub X^+\times X^+$, we write $R^\sharp$ for the congruence on~$X^+$ generated by $R$.  We then define $\pres XR = X^+/R^\sharp$. 
So $\pres XR$ is the set of all $R^\sharp$-classes of words over $X$ under the induced operation on equivalence classes.  We denote the $R^\sharp$-class of $w\in X^+$ by $[w]_R$.  

If $S$ is a semigroup and $\phi:X\to S$ is a function, then $\phi$ extends uniquely to a homomorphism $\al_\phi:X^+\to S$ and an anti-homomorphism $\be_\phi:X^+\to S$ defined, for $x_1\cdots x_k\in X^+$, by
\[
\al_\phi(x_1\cdots x_k) = \phi(x_1)\cdots\phi(x_k) \AND \be_\phi(x_1\cdots x_k) = \phi(x_k)\cdots\phi(x_1).
\]
If $\al_\phi$ preserves $R$, in the sense that $\al_\phi(u)=\al_\phi(v)$ for all $(u,v)\in R$, then $\phi$ extends to a homomorphism
\[
\text{$\ga_\phi:\pres XR\to S$ \qquad defined by \qquad $\ga_\phi([w]_R) = \al_\phi(w)$ \qquad for $w\in X^+$.}
\]
Similarly, if $\be_\phi$ preserves $R$, 
then $\phi$ extends to an anti-homomorphism
\[
\text{$\de_\phi:\pres XR\to S$ \qquad defined by \qquad $\de_\phi([w]_R) = \be_\phi(w)$ \qquad for $w\in X^+$.}
\]
Note that $\al_\phi=\be_\phi$ and (if applicable) $\ga_\phi=\de_\phi$ if $S$ is commutative.

Now let $\Ga$ be a graph with vertex and edge sets $X$ and $E$, respectively.  (Note that $X$ may be an infinite, even uncountable, set.)  The \emph{graph semigroup} on $\Ga$ (see \cite{Fisher1989}) is defined by the presentation
\[
\GS \Ga = \pres X{R_\Ga} \qquad\text{where} \qquad R_\Ga=\set{(xy,yx)}{\{x,y\}\in E}.
\]
As special cases, note that if $\Ga$ has no edges, then $\GS\Ga=X^+$, while if $\Ga=K_X$ is the complete graph on $X$ (which has every possible edge $\{x,y\}$ with $x,y\in X$ and $x\not=y$), then $\GS\Ga$ is the \emph{free commutative semigroup} over $X$, which we denote by~$\FCS X$.  
For simplicity, we will write 
$\wb=[w]_{R_\Ga}$ for the $R_\Ga^\sharp$-class of a word $w\in X^+$.  Note that graph semigroups are also known as \emph{free partially commutative semigroups} or \emph{right angled Artin semigroups} in the literature.
The next result is obvious, and its proof is omitted.

\ms
\begin{lemma}\label{lem:xy=yx}
If $x,y\in X$, then $\overline{xy}=\overline{yx}$ if and only if $\{x,y\}\in E$. \epfres
\end{lemma}

In particular, $\GS\Ga$ is commutative if and only if $\Ga=K_X$, in which case $\GS\Ga=\FCS X$, as noted above.

\ms
\begin{prop}\label{prop:gamma_delta}
Let $\Ga$ be a graph on vertex set $X$, and let $\pi\in\Aut(\Ga)$ be an arbitrary (graph) automorphism of $\Ga$.  Then $\pi$ extends uniquely to a (semigroup) automorphism $\ga_\pi\in\Aut(\GS\Ga)$ and an anti-automorphism $\de_\pi\in\AntiAut(\GS\Ga)$ defined by
\[
\ga_\pi(\overline{x_1\cdots x_k}) = \overline{\pi(x_1) \cdots \pi(x_k)}
\AND
\de_\pi(\overline{x_1\cdots x_k}) = \overline{\pi(x_k) \cdots \pi(x_1)}
\]
for all words $x_1\cdots x_k\in X^+$.
Further, any automorphism or anti-automorphism of $X^+$ is of one of these forms (as appropriate); that is,
\[
\Aut(\GS\Ga) = \set{\ga_\pi}{\pi\in\Aut(\Ga)} \AND
\AntiAut(\GS\Ga) = \set{\de_\pi}{\pi\in\Aut(\Ga)}.
\]
\end{prop}

\pf Again, we just prove the statement for auti-automorphisms.  Let $\pi\in\Aut(\Ga)$.  Since $\{x,y\}\in E$ if and only if $\{\pi(x),\pi(y)\}\in E$, it follows that 
$\de_\pi$ is a well-defined anti-homomorphism, and it is easy to check that $\de_{\pi^{-1}}$ is its inverse mapping.  

Conversely, let ${\ve:\GS\Ga\to\GS\Ga}$ be an arbitrary anti-automorphism.  By Lemma \ref{lem:atoms}, $\ve$ maps $\Atom(\GS\Ga)$ bijectively onto itself, and it is clear that $\Atom(\GS\Ga)= \set{\xb}{x\in X}$.  Define $\pi\in\Sym(X)$ by $\ve(\xb)=\overline{\pi(x)}$.  We claim that, in fact, $\pi\in\Aut(\Ga)$.  Indeed, for any $x,y\in X$, we have
\begin{align*}
\{x,y\}\in E &\iff \overline{xy} = \overline{yx} 
\iff \ve(\overline{xy}) = \ve(\overline{yx}) \\
&\iff \ve(\yb)\ve(\xb) = \ve(\xb)\ve(\yb) 
\iff \overline{\pi(y)\pi(x)} = \overline{\pi(x)\pi(y)} 
\iff \{\pi(x),\pi(y)\}\in E,
\end{align*}
where the first and last step uses Lemma \ref{lem:xy=yx}.  
This completes the proof that $\pi\in\Aut(\Ga)$, and it is clear then that $\ve=\de_\pi$. \epf

\begin{cor}\label{cor:Aut_GSGa}
Let $\Ga$ be a graph on vertex set $X$.  Then
\[
\Aut(\GS\Ga) = \set{\ga_\pi}{\pi\in\Aut(\Ga)} \cong \Aut(\Ga) \AND \AntiAut(\GS\Ga) = \set{\de_\pi}{\pi\in\Aut(\Ga)}.
\]
\begin{itemize}
\itemit{i} If $\Ga=K_X$, then $\AutAntiAut(\GS\Ga)=\Aut(\GS\Ga)\cong\Aut(\Ga) =\Sym(X)$.
\itemit{ii} If $\Ga\not=K_X$, then $\AutAntiAut(\GS\Ga) \cong \Z_2\times\Aut(\Ga)$.
\eit
\end{cor}

\pf By Proposition \ref{prop:gamma_delta}, we have $\Aut(\GS\Ga) = \set{\ga_\pi}{\pi\in\Aut(\Ga)}$.  It is also clear that for any $\pi,\si\in\Aut(\Ga)$, $\ga_\pi=\ga_\si$ if and only if $\pi=\si$, and that $\ga_\pi\ga_\si=\ga_{\pi\si}$, giving $\Aut(\GS\Ga)\cong\Aut(\Ga)$.  If $\Ga=K_X$, then $\GS\Ga=\FCS X$ is commutative, whence $\AutAntiAut(\GS\Ga)=\Aut(\GS\Ga)$, so (i) then follows from the obvious fact that $\Aut(K_X)=\Sym(X)$.  For (ii), note that that for all $\pi,\si\in\Aut(\Ga)$,
\[
\ga_\pi=\ga_\si \iff \de_\pi=\de_\si \iff \pi=\si
\COMMA \qquad
\ga_\pi\ga_\si=\de_\pi\de_\si=\ga_{\pi\si} 
\COMMA \qquad
\ga_\pi\de_\si=\de_\pi\ga_\si=\de_{\pi\si}.
\]
It quickly follows that the map $\Z_2\times\Aut(\Ga) \to \AutAntiAut(\GS\Ga)$ given by $(0,\pi)\mt\ga_\pi$ and $(1,\pi)\mt\de_\pi$ is an isomorphism. \epf

We now wish to use Corollary \ref{cor:Aut_GSGa} to obtain information about the groups $\C(\GS\Ga)=\la I(\GS\Ga)\ra$.  Since $I(\GS\Ga)\sub\AntiAut(\GS\Ga)$, and since $\de_\pi^2=\ga_{\pi^2}$ for $\pi\in\Aut(\Ga)$, we see that
\[
I(\GS\Ga) = \set{\de_\pi}{\pi\in\Aut(\Ga),\ \pi^2=\1}.
\]
Note that $I(\GS\Ga)$ contains $\de_\1$, which is induced by the word-reversing map $X^+\to X^+$.  (So part (ii) of the previous result also follows from Proposition \ref{prop:times}, since $\de_\1$ clearly commutes with any automorphism~$\ga_\pi$.)

For a graph $\Ga$, let us now write $I(\Ga) = \set{\al\in\Aut(\Ga)}{\al^2=\1\not=\al}$ for the set of all \emph{involutions} of $\Ga$, and $\C(\Ga)=\la I(\Ga)\ra$ for the subgroup of $\Aut(\Ga)$ generated by these involutions.  As an immediate consequence of Corollary~\ref{cor:Aut_GSGa} and Lemma \ref{lem:Sym(X)} (which states that $\Sym(X)$ is generated by involutions), we have the following.

\ms
\begin{thm}\label{thm:C(GSGa)}
Let $\Ga$ be a graph on vertex set $X$.
\bit
\itemit{i} If $\Ga=K_X$, then $\C(\GS\Ga)\cong\Sym(X)$.
\itemit{ii} If $\Ga\not=K_X$, then $\C(\GS\Ga)\cong\Z_2\times\C(\Ga)$. \epfres
\eit
\end{thm}

\ms
\begin{rem}
In the special case that $\Ga$ has no edges, we have $\Aut(\Ga)=\Sym(X)$ and, since $\Sym(X)$ is involution generated (Lemma \ref{lem:Sym(X)}), $\C(\Ga)=\Aut(\Ga)$.  We have already noted that $\GS\Ga=X^+$ in this case, so we obtain $\C(X^+)\cong\Z_2\times\Sym(X)$.
\end{rem}

\ms
\begin{rem}
A much simpler construction yields a semigroup $S$ with $\C(S)\cong\Sym(X)$.  We define $S=X\cup\{0\}$, where $0$ is a symbol that does not belong to $X$, and simply declare $uv=0$ for all $u,v\in S$.  Then clearly $\C(S)=\AutAntiAut(S)=\Aut(S)\cong\Sym(X)$.
\end{rem}

As noted in Section \ref{sect:allG}, any group is isomorphic to the automorphism group of a suitable graph.  

\ms
\begin{thm}\label{thm:Z2xG}
Let $G$ be a group generated by involutions.  Then there exists a semigroup $S$ with a proper involution such that $\C(S)\cong\Z_2\times G$.  If $G$ is finite, then~$S$ may be taken to be finitely generated.
\end{thm}

\pf Choose some graph $\Ga$ such that $G\cong\Aut(\Ga)$.  We may assume that $\Ga$ is not complete and is finite if $G$ is finite.  Since $G$ is generated by involutions, so too is $\Aut(\Ga)\cong G$, so $\C(\Ga)=\Aut(\Ga)$.  The result then follows from Theorem \ref{thm:C(GSGa)}(ii). \epf

\section{Other families of semigroups}\label{sect:families}

In this section, we calculate $\C(S)$ for various natural families of semigroups: namely, free groups and free abelian groups of finite rank; finite cyclic groups; symmetric groups; full transformation semigroups; symmetric and dual symmetric inverse monoids; finite diagram monoids; monoids of complex matrices; and rectangular bands.  We also give a natural construction that allows one to embed a semigroup $S$ with $\AntiAut(S)=\emptyset$ into a semigroup $\DS$ with $\AntiAut(\DS)\not=\emptyset$; indeed, the involutions of $\DS$ are in one-one correspondence with the automorphisms of $S$.

\subsection{Free groups of finite rank}\label{subsect:F_X}

Let $X=\{x_1,\ldots,x_n\}$ be a finite set with $|X|=n\geq2$ (the $n=1$ case is considered in the next section), and let $F_X$ be the free group over $X$.  So elements of $F_X$ are equivalence classes of words over $\{x_1,x_1^{-1},\ldots,x_n,x_n^{-1}\}$, where two words are considered equivalent if one may be transformed into the other by successively using the free group relations
\[
x_ix_i^{-1} = x_i^{-1}x_i = 1 \qquad\text{for each $i$,}
\]
where $1$ denotes the empty word.  Nielsen \cite{Nielsen1924} gave a generating set for the automorphism group $\Aut(F_X)$.  Specifically, it was shown that $\Aut(F_X)$ is generated by the automorphisms $\si_1,\ldots,\si_{n-1},\al,\be$ defined by 
\[
\si_i(x_j) = \begin{cases}
x_{i+1} &\text{if $j=i$}\\
x_i &\text{if $j=i+1$}\\
x_j &\text{if $j\not=i,i+1$,}
\end{cases}
\qquad
\al(x_j) = \begin{cases}
x_1^{-1} &\text{if $j=1$}\\
x_j &\text{if $j\not=1$,}
\end{cases}
\qquad
\be(x_j) = \begin{cases}
x_1x_2 &\text{if $j=1$}\\
x_j &\text{if $j\not=1$.}
\end{cases}
\]
(Here we simplify notation by writing a word instead of its equivalence class.)
Actually, in \cite{Nielsen1924}, the automorphisms $\si_2,\ldots,\si_{n-1}$ were replaced by the automorphism $\ga$ defined by $\ga(x_j)=x_{j+1}$ for each~$j$, where the subscript is read modulo $n$.  But since $\{\si_1,\ldots,\si_{n-1}\}$ and $\{\si_1,\ga\}$ both generate the automorphisms of $F_X$ induced by permuting the elements of $X$, the generating set we use here is equivalent.  Note that $\si_1,\ldots,\si_{n-1},\al$ are all of order $2$, but that $\be$ is of infinite order; indeed, $\be^k(x_1)=x_1x_2^k$ for $k\geq1$.  However, it is easy to check that $\be=\be_1\be_2$, where $\be_1,\be_2\in\Aut(F_X)$ are defined by
\[
\be_1(x_j) = \begin{cases}
x_1x_2 &\text{if $j=1$}\\
x_2^{-1} &\text{if $j=2$}\\
x_j &\text{if $j\geq3$,}
\end{cases}
\qquad
\be_2(x_j) = \begin{cases}
x_2^{-1} &\text{if $j=2$}\\
x_j &\text{if $j\not=2$.}
\end{cases}
\]
Moreover, the automorphisms $\be_1,\be_2$ are of order $2$, so it follows that $\Aut(F_X)$ is involution-generated.  Combining this with Proposition \ref{prop:times}, it follows that
\[
\G(F_X)=\la J(F_X)\ra=\Aut(F_X) \AND \C(F_X)= \AutAntiAut(F_X) = \la\iota\ra\times\Aut(F_X),
\]
where $\iota:F_X\to F_X$ is the inversion map, which sends (the equivalence class of) the word $x_{i_1}^{\pm1}\cdots x_{i_k}^{\pm1}$ to $x_{i_k}^{\mp1}\cdots x_{i_1}^{\mp1}$.

\subsection{Free abelian groups of finite rank}\label{subsect:Z^n}

A free abelian group of finite rank $n\geq1$ is the quotient of a free group of rank $n$ by the relations declaring all generators to commute with each other; such a group is isomorphic to $\Z^n$.  Since $\Z^n$ is abelian, $\AutAntiAut(\Z^n)=\Aut(\Z^n)$ and $\C(\Z^n)=\G(\Z^n)$.  The automorphisms of $\Z^n$ are precisely the invertible linear maps $\Z^n\to\Z^n$, where $\Z^n$ is considered as a (free) $\Z$-module.  These are in one-one correspondence with the invertible $n\times n$ matrices over $\Z$: that is, the matrices of determinant $\pm1$.  We denote the set of all such matrices by $\GL(n,\Z)$.  It is well known \cite{Gustafson1991} that $\GL(n,\Z)$ is generated by involutions (i.e., matrices that satisfy $A^2=I_n$, the $n\times n$ identity matrix); Ishibashi \cite{Ishibashi1995} gave an upper bound on the number of involutions required to generate a given element of $\GL(n,\Z)$ for $n\geq3$.  It follows that
\[
\C(\Z^n)=\Aut(\Z^n)=\GL(n,\Z).
\]
It is easy to check that there are infinitely many involutions on $\Z^n$ for $n\geq2$; for example, $\left[\begin{smallmatrix}a&\phantom{-}b\\c&-a\end{smallmatrix}\right]$ is an involution of $\Z^2$ for any $a,b,c\in\Z$ satisfying $a^2+bc=1$.

We note that free groups and free abelian groups have a common generalisation to \emph{graph groups} $\GG\Ga$, defined in an analogous way to the graph semigroups of Section \ref{sect:GS}.  Generators for the automorphism groups $\Aut(\GG\Ga)$ are given in \cite{Laurence1995,Servatius1989}, but a description of the involution-generated subgroup of $\Aut(\GG\Ga)$ is beyond the scope of the present work.

\subsection{Finite cyclic groups}\label{subsect:Z_n}

Let $n\geq2$ be an integer, and 
write $\Z_n=\{0,1,\ldots,n-1\}$ for the additive group of integers modulo $n$.  A homomorphism $\al:\Z_n\to\Z_n$ is uniquely determined by $\al(1)=k$, where $k\in\Z_n$, and any such $k$ gives rise to a homomorphism; so we denote this homomorphism by $\al_k$, and note that $\al_k(m)=km$ for all $m\in\Z_n$.  Since $\al_k\al_l=\al_{kl}$, we see that $\End(\Z_m)$, the semigroup of all endomorphisms of~$\Z_n$ (i.e., all homomorphisms $\Z_n\to\Z_n$), is isomorphic to $(\Z_n,\cdot)$, the multiplicative semigroup of integers modulo~$n$.  The automorphism group $\Aut(\Z_n)$ is therefore isomorphic to $\mathbb U_n=\set{m\in\Z_n}{\gcd(m,n)=1}$, the group of multiplicative units of $(\Z_n,\cdot)$.  The involutions of $\Z_n$ are therefore the maps $\al_k$ with $k^2=1\not=k$ in~$\Z_n$.  Note that the composite of two such involutions is itself an involution or else the identity map.  In other words,
\[
\C(\Z_n) = I(\Z_n)\cup\{\1\} = \set{\al_k}{k\in\Z_n,\ k^2=1},
\]
where as usual $\1\in\Aut(\Z_n)$ denotes the identity automorphism.  Consider the prime factorization $n=2^mp_1^{m_1}\cdots p_r^{m_r}$ (so $m\geq0$, $p_1,\ldots,p_r$ are distinct odd primes, and $m_i\geq1$ for each $i$).  It is well known that the number of solutions to $k^2=1$ in $\Z_n$ is equal to $2^{R(n)}$, where
\[
R(n) = \begin{cases}
r &\text{if $m\leq1$}\\
r+1 &\text{if $m=2$}\\
r+2 &\text{if $m\geq3$.}
\end{cases}
\]
(See for example \cite{Nagell1951}, or \cite{Hage} where involutions on finite abelian groups were studied.)  
Since every element $\al\in\C(\Z_n)=I(\Z_n)\cup\{\1\}$ satisfies $\al^2=\1$, it follows that $\C(\Z_n)\cong\Z_2^{R(n)}$.

\subsection{Symmetric groups}\label{subsect:Sym(X)}

Let $X$ be an arbitrary set.  It is well known \cite{SU1937,Maurer1983} that 
\[
\Aut(\Sym(X)) \cong \begin{cases}
\Sym(X) &\text{if $|X|\not=2,6$} \\
\{\1\}  &\text{if $|X|=2$} \\
\Sym(X)\rtimes\Z_2 &\text{if $|X|=6$.} 
\end{cases}
\]
If $|X|\not=2,6$, then every automorphism of $\Sym(X)$ is inner; that is, if $\al\in\Aut(\Sym(X))$, then there exists some $\pi\in\Sym(X)$ such that $\al(\si)=\pi\si\pi^{-1}$ for all $\si\in\Sym(X)$.  When $|X|=6$, there are non-inner automorphisms sending transpositions to products of three disjoint transpositions, resulting in $\Aut(\Sym(X))\cong\Sym(X)\rtimes\Z_2$; for a nice recent exposition of this fact, see \cite{McCammond2014}.  Since $\Sym(X)$ is a group and is generated by its involutions, Proposition \ref{prop:times} gives $\C(\Sym(X))=\AutAntiAut(\Sym(X))\cong\Z_2\times\Sym(X)$ for $|X|\not=2,6$.  If $|X|=6$, then $\C(\Sym(X))=\AutAntiAut(\Sym(X))\cong\Z_2\times(\Sym(X)\rtimes\Z_2)$.  (If $|X|=2$, then $\C(\Sym(X))=\AutAntiAut(\Sym(X))=\{\1\}$.)

\subsection{Full transformation semigroups}\label{subsect:T_X}

The \emph{full transformation semigroup} on the set $X$, denoted $\T_X$, consists of all transformations of $X$ (i.e., all functions $X\to X$) under the operation of composition.  The symmetric group $\Sym(X)$ is the group of units of $\T_X$.  It is well known \cite{Schreier1937,CW1976,IN1972} that every automorphism of $\T_X$ corresponds to conjugation by some fixed element of $\Sym(X)$, so that $\Aut(\T_X)\cong\Sym(X)$; there are no special cases such as $|X|=6$, here.  However, $\T_X$ has no anti-automorphisms (we assume $|X\geq2|$), as we now explain.  

Recall that \emph{Green's relations} are defined on a semigroup $S$ by
\begin{gather*}
x\R y \iff x\Sone=y\Sone \COMMA x\L y \iff \Sone x=\Sone y \COMMA x\mathscr J y \iff \Sone x\Sone =\Sone y\Sone  , \\
\H=\R\cap\L \COMMA \D=\R\vee\L=\R\circ\L=\L\circ\R.
\end{gather*}
Here, $\Sone$ denotes the monoid obtained from $S$ by adjoining an identity element $1$, if necessary.  It is clear that if $\al:S\to S$ is an anti-automorphism of the semigroup $S$, then for any $x,y\in S$,
\[
x\R y \iff \al(x)\L\al(y) \COMMA 
x\L y \iff \al(x)\R\al(y) \COMMA 
x\mathscr K y \iff \al(x)\mathscr K\al(y) ,
\]
where $\mathscr K$ denotes any of $\H$, $\D$, $\mathscr J$.  In particular, $\al$ maps $\R$-, $\L$- and $\D$-classes of $S$ to $\L$-, $\R$- and $\D$-classes of $S$ (respectively), and if $D$ is a $\D$-class, then the number of $\L$-classes in $D$ is equal to the number of $\R$-classes in $\al(D)$, with a dual statement holding also.

It is well known \cite{Hig,Howie} that for transformations $f,g\in\T_X$,
\begin{gather*}
f\R g \iff \im(f)=\im(g) \COMMA
f\L g \iff \ker(f)=\ker(g) \COMMA
f\D g \iff f\mathscr J g \iff \rank(f)=\rank(g) ,
\end{gather*}
where, as usual, $\im(f)$, $\ker(f)$ and $\rank(f)$ denote the \emph{image}, \emph{kernel} and \emph{rank} of $f\in\T_X$, and are defined by
\[
\im(f) = \set{f(x)}{x\in X} \COMMA
\ker(f) = \set{(x,y)\in X\times X}{f(x)=f(y)} \COMMA
\rank(f) = |\im(f)|.
\]
Note that $\im(f)$ is a subset of $X$, $\ker(f)$ is an equivalence relation on $X$, and $\rank(X)$ is a cardinal (which may be anything from $1$ to $|X|$).  (Note also that it is common in the semigroup literature to write functions to the right of their arguments; with this convention, the roles of the image and kernel in characterising the $\R$ and $\L$ relations would be switched.)  So the $\D=\mathscr J$-classes are given by
\[
D_\lam = \set{f\in\T_X}{\rank(f)=\lam} \qquad\text{for each cardinal $1\leq \lam\leq|X|$.}
\]
  And the $\R$-classes and $\L$-classes contained in the $\D$-class $D_\lam$ are indexed (respectively) by the sets
\[
\mathbb R_\lam=\set{A\sub X}{|A|=\lam} \AND \mathbb L_\lam=\set{\ve}{\text{$\ve$ is an equivalence on $X$ with $\lam$ equivalence classes}},
\]
corresponding (respectively) to the collections of all images and kernels of the given rank $\lam$.
In particular, $|\mathbb R_\lam|\geq2$ and $|\mathbb L_\lam|\geq2$ for any $1\leq\lam\leq|X|$ unless (a) $\lam=1$, where we have $|\mathbb R_1|=|X|$ and $|\mathbb L_1|=1$, or (b) $\lam=|X|<\an$, in which case $|\mathbb R_\lam|=|\mathbb L_\lam|=1$.  So if $\al:\T_X\to\T_X$ was an anti-automorphism, then $\al(D_1)$ would be a $\D$-class of $\T_X$ with just one $\R$-class and $|X|\geq2$ $\L$-classes.  Since there is no such $\D$-class of $\T_X$, it follows that there are no anti-automorphisms of $\T_X$.  That is,
\[
\AntiAut(\T_X)=\emptyset \COMMA
\AutAntiAut(\T_X)=\Aut(\T_X)\cong\Sym(X) \COMMA
\C(\T_X)=\la\emptyset\ra=\{\1\}.
\]
Note that, writing $\FCS X$ for the free commutative semigroup over $X$, as in Section \ref{sect:GS}, we have
\[
\text{$\Aut(\T_X)\cong\Aut(\FCS X)$ \qquad and \qquad $\AutAntiAut(\T_X)\cong\AutAntiAut(\FCS X)$ \qquad but \qquad $\C(\T_X)\not\cong\C(\FCS X)$.}
\]
So the $\C(S)$ invariant enables one to distinguish some semigroups that have isomorphic (signed) automorphism groups but are not isomorphic.  (Of course many other invariants could be used to distinguish between $\T_X$ and $\FCS{X}$.)

\subsection{Symmetric and dual symmetric inverse semigroups}\label{subsect:I_X}

There are natural inverse semigroup analogues of the full transformation semigroup: namely, the \emph{symmetric inverse semigroup} $\I_X$ \cite{Lawson1998,Lipscombe1996}, and the \emph{dual symmetric inverse semigroup} $\I_X^*$ \cite{FL1998}.  These consist, respectively, of all bijections between subsets of $X$ and all bijections between quotients of~$X$, under operations we will not need to describe here.  It is well known \cite{Maltcev2007} that for $|X|\geq3$, $\Aut(\I_X^*)\cong\Sym(X)$, with all automorphisms acting as conjugation by elements of $\Sym(X)\sub\I_X^*$.  
It is also certainly well known that $\Aut(\I_X)\cong\Sym(X)$ for any $X$; though we are unaware of a reference for this, the argument used in \cite{IN1972} to show that $\Aut(\T_X)\cong\Sym(X)$ may easily be adapted to show this.
Since $\I_X$ and $\I_X^*$ are both inverse monoids, the involution $\iota:f\mt f^{-1}$ commutes with every automorphism, so it follows that if $S$ is either $\I_X$ or $\I_X^*$ (with $|X|\geq3$ in the latter case), then
\[
\Aut(S)\cong\Sym(X) \AND \C(S)= \AutAntiAut(S) \cong\Z_2\times\Sym(X).
\]

\subsection{Finite diagram monoids}\label{subsect:P_X}

(Finite) \emph{partition algebras} were introduced independently by Paul Martin \cite{Martin1994} and Vaughan Jones \cite{Jones1994_2} in the context of Potts models in statistical mechanics, and are \emph{twisted semigroup algebras} of the \emph{partition monoids} \cite{Wilcox2007}.  The partition monoid $\P_X$ on the set $X$ contains (isomorphic copies of) the full transformation semigroup $\T_X$ as well as the symmetric and dual symmetric inverse monoids $\I_X$ and~$\I_X^*$ \cite{East2011_2,EF2012,East2011}, and many other diagram monoids \cite{MM2014,FitLau2011,EG2015,LauFit2006,Maz2002,Maltcev2007,Auinger2012,Auinger2014}.  The elements of~$\P_X$ have natural diagrammatic representations, but we will not explicitly describe these (or the algebraic structure of~$\P_X$); for more details, we refer the reader to \cite{EF2012,EG2015,DEEFHHL2015,East2011,Maz2002,ADV2012_2} and references therein.  The partition monoid is a canonical example of a \emph{regular $*$-semigroup} \cite{Nordahl1978}; there is a map ${}^*:\P_X\to\P_X:f\mt f^*$ satisfying $(f^*)^*=f$, $(fg)^*=g^*f^*$ and $f=ff^*f$ (and $f^*=f^*ff^*$) for all $f,g\in\P_X$.  In particular, the~${}^*$ map is an involution; this involution and another were studied in \cite{ADV2012_2} in the context of (inherently) non-finitely based equational theories.  The partition monoid $\P_X$ contains the symmetric group $\Sym(X)$ as its group of units, and if $f\in\Sym(X)$, then $f^*=f^{-1}$.  
It is known \cite{Maz2002} that $\Aut(\P_X)\cong\Sym(X)$ for finite $X$, with all automorphisms corresponding to conjugation by permutations.
Since $(fgf^{-1})^*=f(g^*)f^{-1}$ for all $f\in\Sym(X)$ and $g\in\P_X$, we see that the ${}^*$ map commutes with each automorphism of (finite) $\P_X$, and so
\[
\Aut(\P_X)\cong\Sym(X) \AND \C(\P_X)= \AutAntiAut(\P_X) \cong\Z_2\times\Sym(X).
\]

\subsection{Full linear monoids}\label{subsect:MnC}

Let $\MnC$ denote the monoid of all $n\times n$ matrices over the complex field $\bbC$.  We also write
$\GLnC = \set{A\in\MnC}{\det(A)\not=0}$ for the \emph{general linear group} of degree $n$ over $\bbC$, which is the group of units of $\MnC$.  
For an invertible matrix $A\in\GLnC$, denote by $\phi_A\in\Aut(\MnC)$ the automorphism determined by conjugation by $A$.  That is, $\phi_A:\MnC\to\MnC$ is defined by $\phi_A(X)=AXA^{-1}$ for all $X\in\MnC$.  Putcha \cite{Putcha1983} showed that
\[
\Aut(\MnC)=\set{\phi_A}{A\in\GLnC}.
\]
So $\Phi:\GLnC\to\Aut(\MnC):A\mt\phi_A$ is a surjective group homomorphism, and its kernel, $\ker(\Phi)=\set{A\in\GLnC}{\phi_A=\1}$, is precisely the centre $Z(\GLnC)$ of $\GLnC$; that is, $\ker(\Phi)=Z(\GLnC)=\set{xI_n}{x\in\bbC\sm\{0\}}$.
In particular, $\Aut(\MnC)$ is isomorphic to the quotient $\GLnC/Z(\GLnC)$, which is precisely the \emph{projective general linear group} of degree $n$ over $\bbC$, and denoted $\PGLnC$. 

Denote by $\tau:\MnC\to\MnC:A\mt A^{\operatorname{T}}$ the transpose map.  Since $\tau$ is an anti-automorphism (indeed, an involution), Lemma \ref{lem:commutative} tells us that every anti-automorphism of $\MnC$ is of the form $\phi_A\tau$ for some $A\in\GLnC$; we will write $\psi_A=\phi_A\tau$ for this anti-automorphism.  So
\[
\AutAntiAut(\MnC) 
= \set{\phi_A}{A\in\GLnC}\cup\set{\psi_A}{A\in\GLnC}.
\]
Note that 
\[
\psi_A (X) = AX\tr A^{-1}  \qquad\text{for all $X\in\MnC$.}
\]
One easily checks that $\tau\phi_A=\phi_{(A\tr)^{-1}}\tau=\psi_{(A\tr)^{-1}}$, and that the inverse mapping of $\psi_A$ is $\psi_{A\tr}$.  
In particular,
\[
\psi_A\in I(\MnC) 
\iff
\psi_A=\psi_{A\tr}
\iff
\phi_A=\phi_{A\tr}
\iff
A\tr=xA \quad\text{for some $x\in\bbC\sm\{0\}$}.
\]
So, for example, $\psi_A\in I(\MnC)$ if $A$ is symmetric (i.e., $A=A\tr$).  It is well known \cite{Bosch1986} that any matrix $A\in\MnC$ is the product $A=A_1A_2$ of two symmetric matrices $A_1,A_2\in\MnC$.  In particular, if $A\in\GLnC$, then the symmetric matrices $A_1,A_2$ also belong to $\GLnC$, and we have
\[
\phi_A = \phi_{A_1}\phi_{A_2} =  \phi_{A_1}\tau^2\phi_{A_2} =  \psi_{A_1} \psi_{(A_2\tr)^{-1}} \in \C(\MnC),
\]
since $(A_2\tr)^{-1}=A_2^{-1}$ is also symmetric.  This shows that $\Aut(\MnC)\sub\C(\MnC)$, and it follows that
\[
\C(\MnC) = \AutAntiAut(\MnC) = \Aut(\MnC) \rtimes\la\tau\ra \cong \PGLnC\rtimes\Z_2.
\]
Note here that the semidirect product is not direct, in view of the rule $\tau\phi_A=\phi_{(A\tr)^{-1}}\tau$.  We also note that, since any complex matrix of determinant $1$ is a product of (at most four) involutory matrices \cite{GHR1976}, it quickly follows that $\Aut(\MnC)\cong\PGLnC$ is itself involution-generated; that is, $\Aut(\MnC)=\G(\MnC)=\la J(\MnC)\ra$.

\subsection{Rectangular bands}\label{sect:B}

A \emph{band} is a semigroup consisting entirely of idempotents.  It is well known \cite{Clifford1954} that every band is a semilattice of \emph{rectangular bands}.  A rectangular band is a semigroup $B=X\times Y$ with multiplication $(x_1,y_1)(x_2,y_2)=(x_1,y_2)$ for $x_1,x_2\in X$ and $y_1,y_2\in Y$.  In this section, we fix such a rectangular band $B=X\times Y$, and we describe the group $\C(B)$.

The $\R$-, $\L$- and $\H$-classes of $B$ are precisely the sets
\[
R_x=\{x\}\times Y \COMMA L_y=X\times\{y\} \COMMA H_{xy}=R_x\cap L_y=\{(x,y)\} \qquad\text{for each $x\in X$ and $y\in Y$.}
\]
  Since any anti-automorphism of $B$ maps $\R$-classes to $\L$-classes (and vice versa), $\AntiAut(B)=\emptyset$ if $|X|\not=|Y|$, in which case $\C(B)=\{\id_B\}$.  So we assume from now on that $|X|=|Y|$; in fact, without loss of generality, we may assume that $X=Y$.  So $B=X\times X$ is a \emph{square band}.

For $\si,\tau\in\Sym(X)$, we define mappings 
\[
\ga_{\si,\tau}:B\to B:(x,y)\mt(\si(x),\tau(y)) \AND \de_{\si,\tau}:B\to B:(x,y)\mt(\si(y),\tau(x)).
\]
Note that $\ga_{\1,\1}=\1_B$.  
It is easily seen (and well known in the case of automorphisms) that
\[
\Aut(B)=\set{\ga_{\si,\tau}}{\si,\tau\in\Sym(X)}\cong\Sym(X)\times\Sym(X)  \quad\text{and}\quad \AntiAut(B)=\set{\de_{\si,\tau}}{\si,\tau\in\Sym(X)}.
\]
(More generally, $\Aut(X\times Y)\cong\Sym(X)\times\Sym(Y)$.)  

Put $\iota=\de_{\1,\1}$, so $\iota(x,y)=(y,x)$ for all $x,y\in X$.  Then clearly $\iota\in I(S)$, so Proposition \ref{prop:rtimes} gives
\[
\AutAntiAut(B) = \Aut(B)\rtimes\la\iota\ra \cong (\Sym(X)\times\Sym(X))\rtimes\Z_2 = \Sym(X)\wr\Z_2.
\]
Here $G\wr\Z_2=(G\times G)\rtimes\Z_2$ denotes the \emph{wreath product} in which $\Z_2\cong\Sym(2)$ acts on $G\times G$ by permuting the coordinates.  Specifically, the action of $\la\iota\ra\cong\Z_2$ on $\Aut(B)\cong\Sym(X)\times\Sym(X)$ is given by
\[
\iota\ga_{\si,\tau}\iota = \ga_{\tau,\si} \qquad\text{for all $\si,\tau\in\Sym(X)$.}
\]
We now describe the group $\C(B)=\la I(B)\ra$.  
One easily checks that products in $\AutAntiAut(B)$ are given by 
\[
\ga_{\si_1,\tau_1}\ga_{\si_2,\tau_2} = \ga_{\si_1\si_2,\tau_1\tau_2} ,\ \ \ 
\de_{\si_1,\tau_1}\de_{\si_2,\tau_2} = \ga_{\si_1\tau_2,\tau_1\si_2} ,\ \ \
\ga_{\si_1,\tau_1}\de_{\si_2,\tau_2} = \de_{\si_1\si_2,\tau_1\tau_2} ,\ \ \ 
\de_{\si_1,\tau_1}\ga_{\si_2,\tau_2} = \de_{\si_1\tau_2,\tau_1\si_2} .
\]
In particular, $\de_{\si,\tau}\in I(B)$ is an involution if and only if $\tau=\si^{-1}$, so
\[
I(B) = \set{\de_{\si,\si^{-1}}}{\si\in\Sym(X)}.
\]
By Proposition \ref{prop:rtimes}, $\C(B)=(\C(B)\cap\Aut(B))\rtimes\la\iota\ra$, so we concentrate on the group $\C(B)\cap\Aut(B)$ which, for convenience, we will denote by $G$.  So $G$ consists of all even-length products of elements from $I(B)$.  Note that if $\si,\tau\in\Sym(X)$ are arbitrary, then
\[
\de_{\si,\si^{-1}} \de_{\tau^{-1},\tau} = \ga_{\si\tau,\si^{-1}\tau^{-1}} = \ga_{\si,\si^{-1}} \ga_{\tau,\tau^{-1}},
\]
and that $\ga_{\si,\si^{-1}} = \de_{\si,\si^{-1}}\iota \in G$ (and similarly $\ga_{\tau,\tau^{-1}}\in G$).  So in fact, $G$ is generated by the set 
\[
\set{\ga_{\si,\si^{-1}}}{\si\in\Sym(X)}.
\]
It quickly follows that $G\cong \KSymX$, in the notation of Lemma \ref{lem:KG}.  By the above discussion, and the observations in Example \ref{eg:KSymX}, we see that $\C(B) = \AutAntiAut(B)\cong \Sym(X)\wr\Z_2$ for infinite $X$, and
\[
\C(B) = \set{\ga_{\si,\tau},\de_{\si,\tau}}{\si,\tau\in\Sym(X),\ \si\tau\in\Alt(X)} 
\cong \big((\Alt(X)\times\Alt(X))\rtimes\Z_2\big)\rtimes\Z_2 
\]
for finite $X$.  Let $X$ be finite and put
\[
K = \KSymX = \set{(\si,\tau)\in\Sym(X)\times\Sym(X)}{\si\tau\in\Alt(X)}.
\]
So, as noted in Example \ref{eg:KSymX}, $\Alt(X)\times\Alt(X)\leq K\leq \Sym(X)\times\Sym(X)$.  
In this way, we see that
\[
\Alt(X)\wr\Z_2 \leq K\rtimes\Z_2 \leq \Sym(X)\wr\Z_2,
\]
so $\C(B)\cong K\rtimes\Z_2$ is a \emph{sub-wreath product}; such products are also called \emph{cascade products} \cite{ENN2014}.

\subsection{The doubled semigroup construction}\label{subsect:D_S}

We have seen examples of semigroups with no anti-automorphisms; these include the full transformation semigroups $\T_X$ (Section \ref{subsect:T_X}) and non-square rectangular bands (Section \ref{sect:B}).  In this section, we provide a construction that allows us to extend such a semigroup $S$ with $\AntiAut(S)=\emptyset$ to a semigroup $\DS$ with $\AntiAut(\DS)\not=\emptyset$; indeed, $I(\DS)\not=\emptyset$ and the involutions of $\DS$ are in one-one correspondence with the automorphisms of $S$.

With the above discussion in mind, let $S$ be an arbitrary non-empty semigroup with $\AntiAut(S)=\emptyset$.  (The following construction works also in the case that $\AntiAut(S)\not=\emptyset$, but the calculations of $\AutAntiAut(\DS)$ and $\C(\DS)$ are more involved; the authors can be contacted for details.)  Let $S=\set{s^*}{s\in S}$ be a set disjoint from~$S$, and in one-one correspondence with $S$ via the map $s\mt s^*$.  Define a product on $S^*$ by $s^*t^*=(ts)^*$ for each $s,t\in S$.  So $S^*$ is the \emph{dual semigroup} of $S$.  Now define a new semigroup $\DS=S\cup S^*\cup\{0\}$, where $0$ is a new symbol that does not belong to $S$ (or $S^*$).  We call $\DS$ the \emph{doubled semigroup} obtained from $S$.  Define a product $\cdot$ on $\DS$, for $x,y\in \DS$, by
\[
x\cdot y = \begin{cases}
xy &\text{if $x,y\in S$ or $x,y\in S^*$}\\
0 &\text{otherwise.}
\end{cases}
\]
So $\DS$ contains both $S$ and $S^*$ as subsemigroups.  
Note that $\iota:\DS\to \DS$ defined by
\[
\iota(0)=0 \COMMA \iota(s)=s^* \COMMA \iota(s^*)=s \qquad\text{for all $s\in S$}
\]
is an involution of $\DS$.  In particular, $I(\DS)$ (and hence $\AntiAut(\DS)$) is non-empty.

\ms
\begin{lemma}\label{lem:Aut(T)}
\begin{itemize}
\itemit{i} If $\ga\in\Aut(\DS)$, then $\ga(S)=S$ and $\ga(S^*)=S^*$.
\itemit{ii} If $\de\in\AntiAut(\DS)$, then $\de(S)=S^*$ and $\de(S^*)=S$.
\end{itemize}
\end{lemma}

\pf Since $\iota(S)=S^*$ and $\iota(S^*)=S$, it suffices (by Lemma \ref{lem:commutative}) to prove (i), so suppose $\ga\in\Aut(\DS)$.  First note that if there exists $s,t\in S$ such that $\ga(s)\in S$ and $\ga(t)\in S^*$, then $\ga(0)=0=\ga(s)\ga(t)=\ga(st)$, with $st\not=0$, contradicting the injectivity of $\ga$.  It follows that either (a) $\ga(S)\sub S$ or (b) $\ga(S)\sub S^*$.  A similar argument shows that (a)$^*$ $\ga(S^*)\sub S$ or (b)$^*$ $\ga(S^*)\sub S^*$.  Since $\ga$ is surjective, we cannot have both (a) and (a)$^*$, and neither can we have both (b) and (b)$^*$.  Next suppose (a)$^*$ and (b) hold.  Since $\ga$ is surjective, we would have $\ga(S)=S^*$ (and $\ga(S^*)=S$).  But then the restriction of $\ga$ to~$S$ would yield an isomorphism $S\to S^*$, and composing this with the anti-isomorphism $S^*\to S:s^*\mt s$ would yield an anti-automorphism of $S$, contradicting the assumption that $\AntiAut(S)=\emptyset$.  It follows that (a) and (b)$^*$ hold.  Again, the surjectivity of $\ga$ shows that $\ga(S)=S$ and $\ga(S^*)=S^*$. \epf

For automorphisms $\al,\be\in\Aut(S)$, we define maps $\ga_{\al,\be},\de_{\al,\be}:\DS\to\DS$ by
\[
\begin{array}{rclrclrcll}
\ga_{\al,\be}(0)&\!\!\!=\!\!\!&0\COMMA &\ga_{\al,\be}(s)&\!\!\!=\!\!\!&\al(s) \COMMA &\ga_{\al,\be}(s^*)&\!\!\!=\!\!\!&\be(s)^* &\qquad\text{for all $s\in S$}, \\
\de_{\al,\be}(0)&\!\!\!=\!\!\!&0\COMMA &\de_{\al,\be}(s)&\!\!\!=\!\!\!&\be(s)^* \COMMA &\de_{\al,\be}(s^*)&\!\!\!=\!\!\!&\al(s) &\qquad\text{for all $s\in S$}.
\end{array}
\]
It is easy to check that $\ga_{\al,\be}\in\Aut(\DS)$ and $\de_{\al,\be}\in\AntiAut(\DS)$ for all $\al,\be\in\Aut(S)$.  It then follows quickly from Lemma \ref{lem:Aut(T)} that
\[
\Aut(\DS)=\set{\ga_{\al,\be}}{\al,\be\in\Aut(S)} \AND \AntiAut(\DS) = \set{\de_{\al,\be}}{\al,\be\in\Aut(S)}.
\]
It is easily checked that for any $\al_1,\al_2,\be_1,\be_2\in\Aut(S)$,
\[
\ga_{\al_1,\be_1}\ga_{\al_2,\be_2} = \ga_{\al_1\al_2,\be_1\be_2} ,\ \ 
\de_{\al_1,\be_1}\de_{\al_2,\be_2} = \ga_{\al_1\be_2,\be_1\al_2} ,\ \
\ga_{\al_1,\be_1}\de_{\al_2,\be_2} = \de_{\al_1\al_2,\be_1\be_2} ,\ \ 
\de_{\al_1,\be_1}\ga_{\al_2,\be_2} = \de_{\al_1\be_2,\be_1\al_2} .
\]
Since also
$
\iota\ga_{\al,\be}\iota=\ga_{\be,\al} 
$
for all $\al,\be\in\Aut(S)$, it follows that
\[
\Aut(\DS)\cong\Aut(S)\times\Aut(S) \AND \AutAntiAut(\DS) = \Aut(\DS)\rtimes\la\iota\ra \cong 
\Aut(S)\wr\Z_2.
\]
In a similar fashion to Section \ref{sect:B}, we have
\[
I(\DS) = \set{\de_{\al,\al^{-1}}}{\al\in\Aut(S)},
\]
and it again quickly follows that
\[
\C(\DS)\cap\Aut(\DS) \cong \KAutS \AND \C(\DS) \cong \KAutS\rtimes\Z_2.
\]
Again, we note that 
$
\der{\Aut(S)}\wr\Z_2 \leq \KAutS\rtimes\Z_2 \leq \Aut(S)\wr\Z_2,
$
so that $\C(\DS)\cong \KAutS\rtimes\Z_2$ is again a cascade product \cite{ENN2014}.

\ms
\begin{eg}
Consider the case in which $S=\T_X$ is the full transformation semigroup on the set $X$.  We saw in Section \ref{subsect:T_X} that $\AntiAut(\T_X)=\emptyset$ and $\Aut(\T_X)\cong\Sym(X)$.  It then follows that
\[
\Aut(\DTX) \cong \Sym(X)\times\Sym(X)  \COMMA \AutAntiAut(\DTX) \cong \Sym(X)\wr\Z_2 \COMMA \C(\DTX)\cong \KSymX\rtimes\Z_2.
\]
The group $\KSymX$ is described in Example~\ref{eg:KSymX}.  Interestingly, by consulting Section \ref{sect:B}, we see that $\mathcal A(\DTX)\cong\mathcal A(B)$, where $B=X\times X$ is a square band and $\mathcal A$ denotes any of the operators $\Aut$, $\AutAntiAut$,~$\C$.  We obtain the same isomorphisms in the case that $S=X$ is a left-zero semigroup (in which case, $S^*=X^*$ is a right zero semigroup).
\end{eg}

\subsection{Summary}

Table \ref{tab:C(S)} summarises the calculations of $\C(S)$ for the various semigroups $S$ we have considered in Sections~\ref{sect:GS} and \ref{sect:families}.

\setlength{\tabcolsep}{1em}
\begin{table}[H]%
\scalebox{0.7}{\parbox{\linewidth}{%
\begin{center}
\begin{tabular}{l|rrrr}
Semigroup $S$ & $\Aut(S)$ & $\AutAntiAut(S)$ & $\C(S)$ & Section\\
\hline
Free  semigroup, $X^+$ & $\Sym(X)$ & $\Z_2\times\Sym(X)$ & $\Z_2\times\Sym(X)$ & \ref{sect:GS} \\
Free commutative semigroup, $\FCS X$ & $\Sym(X)$ & $\Sym(X)$ & $\Sym(X)$ & \ref{sect:GS} \\
Graph semigroup, $\GS\Ga$ ($\Ga\not=K_X$) & $\Aut(\Ga)$ & $\Z_2\times\Aut(\Ga)$ & $\Z_2\times\C(\Ga)$ & \ref{sect:GS} \\
Free group, $F_X$, $|X|<\an$ & $\Aut(F_X)$ & $\Z_2\times\Aut(F_X)$ & $\Z_2\times\Aut(F_X)$ & \ref{subsect:F_X} \\
Free abelian group, $\Z^n$ & $\GL(n,\Z)$ & $\GL(n,\Z)$ & $\GL(n,\Z)$ & \ref{subsect:Z^n} \\
Finite cyclic group, $\Z_n$ & $\mathbb U_n$ & $\mathbb U_n$ & $\Z_2^{R(n)}$ & \ref{subsect:Z_n} \\
Symmetric group, $\Sym(X)$, $|X|\not=1,2,6$ & $\Sym(X)$ & $\Z_2\times\Sym(X)$ & $\Z_2\times\Sym(X)$ & \ref{subsect:Sym(X)} \\
Symmetric group, $\Sym(X)$, $|X|=6$ & $\Sym(X)\rtimes\Z_2$ & $\Z_2\times(\Sym(X)\rtimes\Z_2)$ & $\Z_2\times(\Sym(X)\rtimes\Z_2)$ & \ref{subsect:Sym(X)} \\
Full transformation semigroup, $\T_X$ & $\Sym(X)$ & $\Sym(X)$ & $\{\1\}$ & \ref{subsect:T_X} \\
Symmetric inverse monoid, $\I_X$ & $\Sym(X)$ & $\Z_2\times\Sym(X)$ & $\Z_2\times\Sym(X)$ & \ref{subsect:I_X} \\
Dual symmetric inverse monoid, $\I_X^*$, $|X|\geq3$ & $\Sym(X)$ & $\Z_2\times\Sym(X)$ & $\Z_2\times\Sym(X)$ & \ref{subsect:I_X} \\
Partition monoid, $\P_X$, $|X|<\an$ & $\Sym(X)$ & $\Z_2\times\Sym(X)$ & $\Z_2\times\Sym(X)$ & \ref{subsect:P_X} \\
Full linear monoid, $\MnC$ & $\PGLnC$ & $\PGLnC\rtimes\Z_2$ & $\PGLnC\rtimes\Z_2$ & \ref{subsect:MnC} \\
Rectangular band $X\times Y$, $|X|\not=|Y|$ & $\Sym(X)\times\Sym(Y)$ & $\Sym(X)\times\Sym(Y)$ & $\{\id\}$ & \ref{sect:B} \\
Square band $X\times X$, $|X|\geq\an$ & $\Sym(X)\times\Sym(X)$ & $\Sym(X)\wr\Z_2$ & $\Sym(X)\wr\Z_2$ & \ref{sect:B} \\
Square band $X\times X$, $|X|<\an$ & $\Sym(X)\times\Sym(X)$ & $\Sym(X)\wr\Z_2$ & $((\Alt(X)\times\Alt(X))\rtimes\Z_2)\rtimes\Z_2$ & \ref{sect:B} \\
Doubled semigroup $\DS$, $\AntiAut(S)=\emptyset$ & $\Aut(S)\times\Aut(S)$ & $\Aut(S)\wr\Z_2$ & $\KAutS\rtimes\Z_2$ & \ref{subsect:D_S} \\
Doubled semigroup $\DTX$, $|X|\geq\an$ & $\Sym(X)\times\Sym(X)$ & $\Sym(X)\wr\Z_2$ & $\Sym(X)\wr\Z_2$ & \ref{subsect:D_S} \\
Doubled semigroup $\DTX$, $|X|<\an$ & $\Sym(X)\times\Sym(X)$ & $\Sym(X)\wr\Z_2$ & $((\Alt(X)\times\Alt(X))\rtimes\Z_2)\rtimes\Z_2$ & \ref{subsect:D_S} \\
\end{tabular}
\end{center}
}}
\caption{Semigroups $S$, along with $\Aut(S)$, $\AutAntiAut(S)$ and $\C(S)$.  See the specified sections for further information.}
\label{tab:C(S)}
\end{table}

\section{Open problems}

We conclude the article with a short list of open problems that we believe are worth investigating.  

\ms
\begin{problem}
Classify the (involution-generated) groups $G$ for which $G\cong\C(S)$ for some semigroup $S$ with a proper involution.
\end{problem}

\ms
\begin{problem}
Classify the (involution-generated) groups $G$ for which $G\cong\C(H)$ for some group $H$.  
\end{problem}

\ms
\begin{problem}
Classify the semigroups $S$ for which $\C(S)=\AutAntiAut(S)$. 
\end{problem}

\ms
\begin{problem}
Describe the group $\C(S)$ in the case that $S$ is one of the following: a band; a completely $0$-simple semigroup; a semigroup of matrices over any field; a finite abelian group; a free group or free abelian group of infinite rank; a graph group; an infinite partition monoid.
\end{problem}

\footnotesize\def\bibspacing{-1.1pt}
\bibliography{involutions_biblio}
\bibliographystyle{plain}\end{document}